\documentclass[letterpaper, 10 pt, conference]{ieeeconf}  

\usepackage{float}
\usepackage{cite}
\usepackage{amsmath, amssymb, mathtools}
\usepackage{algorithmic}
\usepackage{comment}
\usepackage{graphicx}
\usepackage{color}

\newtheorem{theorem}{Theorem}
\newtheorem{lemma}{Lemma}%

\newenvironment{pf}
{\noindent\textit{Proof.}\hspace{1.2mm}}
{\hspace{\fill}$\square$}
\DeclareMathOperator{\ber}{ber}
\DeclareMathOperator{\bei}{bei}

\IEEEoverridecommandlockouts                              

\title{\LARGE \bf
Extremum-Seeking Boundary Control for Schrödinger-Type PDEs*
}

\author{Paulo Henrique F. Biazetto$^{1}$, Gustavo Artur de Andrade$^{1}$, Tiago Roux Oliveira$^{2}$ and Miroslav Krstic$^{3}$
\thanks{This work was partially supported by CAPES under grants 88887.629803/2021-00  and 88881.878833/2023-01 (SticAmSud).}
\thanks{$^{1}$Paulo Henrique F. Biazetto and Gustavo A. de Andrade
are with the Department of Automation and Systems, Federal University of Santa Catarina (UFSC), Florianópolis, SC 88040-900, Brazil {\tt\small paulo.biazetto@posgrad.ufsc.br, gustavo.artur@ufsc.br}}%
\thanks{$^{2}$Tiago Roux Oliveira is with the Department of Electronics and Telecommunication Engineering, State University of Rio de Janeiro (UERJ), RJ 20550-900, Brazil {\tt\small tiagoroux@uerj.br}}%
\thanks{$^{3}$Miroslav Krstic is with the Department of Mechanical and Aerospace Engineering, University of California (UCSD), San Diego, CA 92093-0411, USA {\tt\small krstic@ucsd.edu}}
}

\begin{document}

\maketitle
\thispagestyle{empty}
\pagestyle{empty}

\begin{abstract}

This paper addresses the design and analysis of an extremum-seeking (ES) controller for scalar static maps in the context of infinite-dimensional dynamics governed by complex-valued partial differential equations (PDEs) of Schrödinger type. The system is actuated at one boundary, and the map input is defined as a real-valued quadratic functional corresponding to the squared norm of the complex state at the uncontrolled boundary. An isomorphism between the complex Hilbert space and its two-dimensional real-valued representation is established to enable the use of the standard multivariable Newton-based ES method. To compensate for the PDE actuation dynamics, a boundary control strategy based on a two-step backstepping procedure is employed. With a perturbation-based estimate of the Hessian’s inverse, the local exponential stability to a small neighborhood of the unknown extremum point is proved. A numerical example illustrates the effectiveness of the proposed extremum-seeking methodology.
\end{abstract}

\section{Introduction}

The boundary control of the Schrödinger equation has attracted considerable attention within the control theory of distributed parameter systems due to its relevance in quantum dynamics, wave propagation, and flexible structures governed by dispersive partial differential equations (PDEs). Several boundary control strategies have been proposed to achieve stabilization, regulation, and observation of such systems. Early works focused on exact controllability via multiplier techniques and spectral analysis~\cite{phung2001}, whereas more recent approaches have emphasized feedback stabilization by exploiting the dissipative mechanisms introduced through boundary conditions~\cite{steeves2020,Ren2013,Bhandari2024}. Dynamic boundary feedback laws—often involving interconnections with auxiliary systems such as heat or wave equations—have been shown to induce exponential or Gevrey-type stability of the closed-loop dynamics~\cite{Wang2012}. These designs rely on the spectral properties of the underlying Schrödinger operator, whose eigenvalues are typically located on the imaginary axis. By introducing suitable feedback, the spectrum can be shifted toward the left half-plane, ensuring asymptotic or exponential stability.

Beyond stabilization, the optimal control of Schrödinger-type equations has emerged as a key topic at the interface of control theory and quantum physics. Applications span quantum optics, material science, and quantum computing, where the goal is to manipulate quantum states by determining optimal external fields or potential profiles that drive the system toward a desired configuration~\cite{Yildiz2001,Breckner2021}. Within this context, reference \cite{AronnaBonnansKroner2016} formulated optimality conditions for control problems in complex Hilbert spaces and demonstrated their applicability to bilinear Schrödinger dynamics, while publication \cite{Ito2007} provided foundational results on bilinear optimal control in abstract Schrödinger settings. These ideas have been extended to stochastic and nonlinear cases~\cite{CuiLiuZhou2023,Wang2018}, as well as to small-time controllability of bilinear quantum systems~\cite{Chambrion2023}. Parallel advances in quantum technologies have introduced high-fidelity optimization techniques for finite-dimensional approximations of Schrödinger dynamics~\cite{Koch2022,DeKeijzer2025}, enriching the theoretical foundation for the analysis and control of dispersive PDEs.

Motivated by these developments, the present work proposes an extremum seeking control (ESC) scheme for a linearized Schrödinger equation under boundary actuation. This is the first instance of ESC being applied to complex-valued PDEs. The objective is to steer the system toward an optimal steady state corresponding to the extremum of an unknown performance map, while compensating for the PDEs' inherent phase shift and dispersive dynamics through a boundary-based feedback structure. The proposed approach extends the ESC framework introduced in~\cite{book2022}, where the compensator for the known actuation dynamics is designed via an infinite-dimensional backstepping transformation and driven by online estimates of the unknown gradient and Hessian of the static map. In the present setting, the static map is defined as a quadratic function of the real and imaginary parts of the system state. By exploiting the isomorphism between the complex plane and $\mathbb{R}^2$, the real-time optimization problem is reformulated in the real plane, allowing the application of a Newton-based tuning law for Hessian's inverse estimation \cite{Ghaffari2011MultivariableNE}. Invoking the averaging theorem for infinite-dimensional systems \cite{Hale1990}, we establish the asymptotic convergence of the closed-loop dynamics to a neighborhood of the extremum point, while preserving exponential stability of the underlying Schrödinger system. Numerical simulations are provided to illustrate the effectiveness of the proposed method and confirm its convergence properties.

The paper is organized as follows. Section \ref{section:problem} introduces the linearized Schrödinger PDE and the corresponding control objectives with ES. In Section \ref{section:extremum_seeking_design}, we present the proposed ESC design. We begin by designing the demodulation and additive probing signals. Next, we derive the error dynamics and design a compensator using a backstepping methodology. The closed-loop stability and asymptotic convergence to the extremum are analyzed in Section \ref{section:stability}. Section \ref{section:results} illustrates the proposed control system through a numerical simulation example. Finally, Section \ref{section:conclusion} brings the concluding remarks. 

\section{Problem Formulation}\label{section:problem}
\subsection{Schrodinger Mathematical Model}
In this work, we consider the following linearized Schrodinger equation:
\begin{align}
    &u_t(t,x) + i\,u_{xx}(t,x)=0, \label{eq:schrodingerPDE}\\
    &u(t,1) = \theta(t), \label{eq:control_schrodinger}\\
    &u_{x}(t,0) = 0,\label{eq:schrodinger_bc0}\\
    &u(0,x) = u_0 (x),\label{eq:schrodinger_in}
\end{align}
\noindent where $t\in [0,\, +\infty)$ is the time, $x\in[0,\, 1]$ is the space, $u$ is the complex-valuated state, $i$ is the imaginary unit, and $\theta\in \mathbb{C}$ is the control input. The initial condition is given by $u_0 (x)\in H^{2}(0,1)$.

\subsection{Control Problem}\label{section:control_problem}

The goal of the ES method is to minimize an unknown static map \( y = Q(\Theta) \) through real-time optimization, where \( y^{\star} \) and \( \Theta^{\star} \) denote the unknown maximum (minimum) and its maximizer (minimizer), respectively. It is assumed that the map is modeled as a real-valued quadratic functional on the complex Hilbert space \( (\mathbb{C}, \langle \cdot, \cdot \rangle_{\mathbb{R}}) \), with  $\langle z_1, z_2 \rangle_{\mathbb{R}} = \Re(\bar{z}_1z_2)$, for all $z_{1},z_{2}\in\mathbb{C}$, and 
\begin{align}
    \Theta(t) = u(t,0).
\end{align}

Locally around \( \Theta^{\star} \), the map is assumed quadratic:
\begin{equation}
Q(\Theta(t)) = y^\star + \tfrac{1}{2}\langle \Theta(t) - \Theta^\star,\ \mathcal{H}(\Theta(t) - \Theta^\star) \rangle_{\mathbb{R}},
\end{equation}
where \( \mathcal{H} \) is a bounded, self-adjoint, and negative definite operator (or positive definite, for minimization).  

Using the isomorphism \( \psi(\Theta) = (\Re\{\Theta\},\,\Im\{\Theta\}) \), with inverse $\psi^{-1}(h)=h_{1}+ih_{2}$, with $h=(h_{1},h_{2})$ and $h_{1},h_{2}\in\mathbb{R}$,  we obtain the equivalent real representation
\begin{equation}
Q(\Theta(t)) = y^\star + \tfrac{1}{2}(\psi(\Theta(t)) - \psi(\Theta)^\star)^\top H(\psi(\Theta(t)) - \psi(\Theta^\star)), 
\label{extrachapter.eq:static_map}
\end{equation}
where \( H \in \mathbb{R}^{2\times 2} \) is symmetric and negative definite matrix for maximization problems (or positive definite, for minimization).  
This formulation treats the real and imaginary parts of \( \Theta = (\Re\{\Theta\},\Im\{\Theta\}) \) as independent variables, enabling the ES controller design in \( \mathbb{R}^2 \) while preserving consistency with the original complex PDE dynamics.

\section{Extremum Seeking Boundary Control Design} \label{section:extremum_seeking_design}
\subsection{Demodulation Signals}

The demodulation signal $N$, which is used to estimate the Hessian of the static map by multiplying it with the output $y(t)$ of the static map, is defined as
\begin{align}
\hat{H}(t) = N(t)y(t),
    \label{eq:hessian}
\end{align}
\noindent where the elements of the $2\times 2$ demodulating matrix $N$ are
\begin{align*}
    N_{11}(t) &= \frac{16}{a_{1}^2}\left( \sin (\omega_{1} t) - \frac{1}{2} \right),\\
    N_{12}(t) &= N_{21}(t) = \frac{4}{a_{1}a_{2}} \sin (\omega_{1} t) \sin (\omega_{2} t),\\
    N_{22}(t) &= \frac{16}{a_{2}^2}\left( \sin (\omega_{2} t) - \frac{1}{2} \right),
\end{align*}
\noindent with $\omega_{1} \neq \omega_{2}$ and $a_1$ and $a_2$ positive constants. The probing frequencies $\omega_{i}$'s can be selected as $\omega_{i}=\omega'_{i}\omega=\mathcal{O}(\omega)$, for $i\in\{1,\,2\}$, where $\omega$ is a positive constant and $\omega_{i}'$ is a rational number. 

The estimation of the Hessian matrix in non-model-based optimization represents a key challenge, as the Newton-type control law requires not only an accurate estimate of the Hessian but also its inverse. Since the estimated matrix evolves in continuous time and may not be invertible at every instant, a dynamic approach is adopted to generate the inverse asymptotically. Following a Riccati-type formulation, the Hessian inverse estimator is governed by
\begin{equation}
\dot{\Gamma}(t) = \omega_r \Gamma(t) - \omega_r \Gamma(t) \hat{H}(t) \Gamma(t),
\label{eq:gamma_dot}
\end{equation}
where $\omega_r > 0$ is a design constant. 

Defining the estimation error of the inverse as 
\begin{equation}
\tilde{\Gamma}(t) = \Gamma(t) - H^{-1},
\label{eq:gamma_tilde}
\end{equation}
its dynamics can be derived as
\begin{equation}
    \dot{\tilde{\Gamma}}(t) = \omega_r (\tilde{\Gamma}(t) + H^{-1})[I_{2\times 2} - \hat{H}(t)(\tilde{\Gamma}(t) + H^{-1})].
    \label{eq:gamma_tilde_dot}
\end{equation}

This formulation guarantees that $\Gamma(t)$ converges asymptotically to $H^{-1}$, while allowing the convergence rate to be tuned via $\omega_r$. The measurable signal
\begin{equation}
z(t) = \Gamma(t) G(t),
\label{eq:z_def}
\end{equation}
is introduced to derive the control law in the subsequent analysis.

In \eqref{eq:z_def}, $G$ represents the gradient estimate and is given as \cite{book2022}
\begin{align}
    G(t) &= M(t)y(t),\nonumber\\
    M(t) &= \left(\frac{2}{a_{1}}\sin(\omega_{1} t),\;\; \frac{2}{a_{2}}\sin(\omega_{2} t) \right). \label{eq:gradient_estimate}
\end{align}

\subsection{Additive Probing Signal}
The perturbation signal $S$ is adapted from the standard ES framework to the case of PDE actuation dynamics \cite{Oliveira2017ExtremumSF,book2022}. The trajectory generation problem, following \cite{Krsti2008BoundaryCO}, is formulated as
\begin{align}
    R_{t}(t,x) + i\,R_{xx}(t,x) &= 0, \label{eq:trajectoryPDE}\\
    R_{x}(t,0) &= 0, \label{eq:trajectorybc0}\\
    R(t,0) &= a_{1}\sin(\omega_{1} t) +i\,a_{2}\sin(\omega_{2} t), \label{eq:trajectorybc03}\\
    S(t) &\coloneqq R(t,1). \label{eq:perturbation}
\end{align}

The explicit solution of \eqref{eq:trajectoryPDE}--\eqref{eq:perturbation} is established in the following lemma.

\begin{lemma}
    The solution of problem \eqref{eq:trajectoryPDE}--\eqref{eq:trajectorybc03} is 
    \begin{align*}
        R(t,x) = R_r(t,x) + iR_{i}(t,x),
    \end{align*}
    \noindent where
\begin{align*}
    R_r(t,x) &= \frac{a_{1}}{2}\sin (\omega_{1}t)\left[\cos (\sqrt{\omega_{1}}\, x)+\cosh(\sqrt{\omega_{1}}\, x)\right] \\
    & + \frac{a_{2}}{2}\cos (\omega_{2}t)\left[\cos (\sqrt{\omega_{2}}\, x) - \cosh(\sqrt{\omega_{2}}\, x)\right],\\
    R_i(t,x) &= -\frac{a_{1}}{2}\cos (\omega_{1}t)\left[\cos (\sqrt{\omega_{1}}\, x)-\cosh(\sqrt{\omega_{1}}\, x)\right] \\
    & + \frac{a_{2}}{2}\sin (\omega_{2}t)\left[\cos (\sqrt{\omega_{2}}\, x)+\cosh(\sqrt{\omega_{2}}\, x)\right].
\end{align*}
\end{lemma}
\begin{pf}
Let $R(t,x)$ be the solution of \eqref{eq:trajectoryPDE}--\eqref{eq:trajectorybc03}, and consider the power series expansion
\begin{align}
    R(t,x) = \sum_{k=0}^{\infty} a_{k}(t)\frac{x^{k}}{k!}.
    \label{eq:series_trajectory_postulation}
\end{align}

Substituting \eqref{eq:series_trajectory_postulation} into \eqref{eq:trajectoryPDE}--\eqref{eq:trajectorybc03} yields the recurrence relation $a_{k+2}(t) = i\,\dot{a}_{k}(t)$.

From the boundary conditions \eqref{eq:trajectorybc0}--\eqref{eq:trajectorybc03}, one obtains
\begin{equation*}
    a_{0}(t) = a_{1}\sin(\omega_{1}t) + i\,a_{2}\sin(\omega_{2}t),
    \qquad 
    a_{1}(t) = 0.
\end{equation*}

Consequently, by the recurrence relation, all the odd coefficients vanish, i.e., $a_{2k+1}(t) = 0$,  $k \ge 0$.

For the even indices, repeated use of the recurrence formula gives $a_{2k}(t) = i^{\,k}\,a_{0}^{(k)}(t)$, and therefore
$R(t,x) = a_{1}\,S_{\omega_{1}}(t,x) + i\,a_{2}\,S_{\omega_{2}}(t,x)$,
where
\begin{equation*}
    S_{\omega}(t,x) = \sum_{k=0}^{\infty} i^{\,k}\,\sin^{(k)}(\omega t)\,\frac{x^{2k}}{(2k)!}.
\end{equation*}

Using the derivatives of the sine function, this expression can be rewritten as
\begin{align*}
    S_{\omega}(t,x) &= \sin(\omega t)\sum_{m=0}^{\infty} \frac{(\omega x^{2})^{2m}}{(4m)!} \\
    & + i\,\omega\cos(\omega t)\sum_{m=0}^{\infty} \frac{(\omega x^{2})^{2m+1}}{(4m+2)!}.
\end{align*}

The two series above correspond to the well-known expansions
\begin{align*}
    \sum_{m=0}^{\infty}\frac{(\omega x^{2})^{2m}}{(4m)!}
    &= \tfrac{1}{2}\left[\cos(\sqrt{\omega}\,x) + \cosh(\sqrt{\omega}\,x)\right],\\
    \sum_{m=0}^{\infty}\frac{(\omega x^{2})^{2m+1}}{(4m+2)!}
    &= \tfrac{1}{2\sqrt{\omega}\,x}\left[\cosh(\sqrt{\omega}\,x) - \cos(\sqrt{\omega}\,x)\right].
\end{align*}

Substituting these identities into the expression for $S_{\omega}(t,x)$ and then into $R(t,x)$ yields the desired closed-form expressions for $R_{r}(t,x)$ and $R_{i}(t,x)$. Hence, the result follows.
\end{pf}

\subsection{Estimation Errors and PDE-Error Dynamics}

Since our goal is to determine $\Theta^{*}$, which corresponds to the optimal but unknown actuator $\theta(t)$, we introduce the following estimation errors:
\begin{align}
    \hat{\theta}(t) &\coloneqq \theta(t) - S(t), \label{eq:error_sigma}\\
    \hat{\Theta}(t) &\coloneqq \Theta(t) - (a_1\sin(\omega_1 t)+i\,a_2\sin(\omega_2 t)). \label{eq:error_theta1} 
\end{align}

\begin{figure}[H]
\vspace{5pt}
    \includegraphics[width=\linewidth]{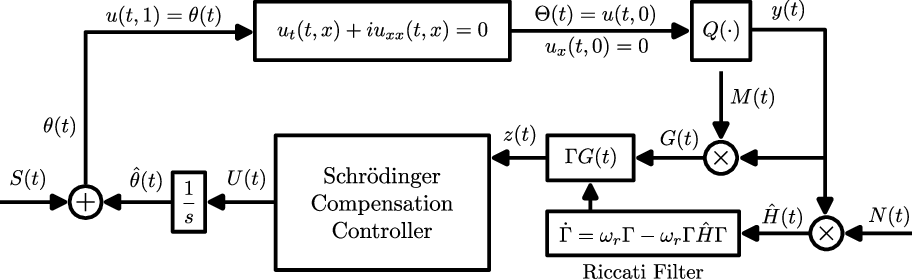}
    \caption{Block diagram of the ESC system applied to the Schrödinger problem.}\label{fig:schrodinger}
\end{figure}

We further define the estimation errors in both the input and the propagated input variables as
\begin{align}
    \tilde{\theta}(t) &\coloneqq \hat{\theta}(t) - \Theta^{*}, \label{eq:error_thetatilde}\\
    \vartheta(t) &\coloneqq \hat{\Theta}(t) - \Theta^{*}. \label{eq:error_vartheta}
\end{align}

Next, consider the state error
\begin{equation}
    \alpha(t,x) \coloneqq u(t,x) - R(t,x) - \Theta^{*}. \label{eq:errorstate}
\end{equation}

Differentiating \eqref{eq:errorstate} with respect to time and substituting \eqref{eq:schrodingerPDE} and \eqref{eq:trajectoryPDE}, we obtain the error dynamics
\begin{align}
    \alpha_{t}(t,x) + i\,\alpha_{xx}(t,x) = 0. \label{eq:errorPDE}
\end{align}

Differentiating \eqref{eq:errorstate} with respect to space, evaluating at $x=0$, and applying the boundary conditions \eqref{eq:schrodinger_bc0} and \eqref{eq:trajectorybc0}, yields
\begin{align}
    \alpha_{x}(t,0) = 0.
\end{align}

From \eqref{eq:error_vartheta} and the definition in \eqref{eq:errorstate}, we also obtain
\begin{align}
    \vartheta(t) = \alpha(t,0). \label{eq:erroTheta}
\end{align}

Similarly, evaluating \eqref{eq:errorstate} at $x=1$, and using the boundary condition \eqref{eq:control_schrodinger} together with \eqref{eq:perturbation}, \eqref{eq:error_sigma} and \eqref{eq:error_thetatilde}, we have
\begin{align}
    \alpha(t,1) = \tilde{\theta}(t). \label{eq:errorbc1}
\end{align}

Taking the time derivative of \eqref{eq:errorPDE}--\eqref{eq:errorbc1} and defining $U(t)\coloneqq \dot{\tilde{\theta}}(t)$, the propagated error dynamics is given by
\begin{align}
    \dot{\vartheta}(t) &= \beta(t,0), \label{eq:ODE1errordynamics}\\
    \beta_{t}(t,x) + i\beta_{xx}(t,x) &= 0, \label{eq:PDEerrordynamics}\\
    \beta_{x}(t,0) &= 0, \label{eq:BC0errordynamics}\\
    \beta(t,1) &= U(t), \label{eq:BCUerrordynamics}
\end{align}
where $\beta(t,x) \coloneqq \alpha_{t}(t,x)$.

Adapting the scheme proposed in \cite{book2022} and combining \eqref{eq:schrodingerPDE}--\eqref{eq:schrodinger_in} and \eqref{extrachapter.eq:static_map}, the closed-loop extremum seeking controller with actuation dynamics governed by the Schrödinger PDE is illustrated in Figure~\ref{fig:schrodinger}.

\subsection{Schrodinger Compensation via Backstepping Boundary Control}

\subsubsection{Target System}

We aim to map the system \eqref{eq:ODE1errordynamics}--\eqref{eq:BCUerrordynamics} into the exponentially stable ODE--PDE cascade
\begin{align}
    \dot{\vartheta}(t) &= -K\vartheta(t) + w(t,0), \label{eq:target_ode}\\
    w_{t}(t,x) +i\,w_{xx}(t,x)&= - c w(t,x), \label{eq:target_pde}\\
    w_{x}(t,0) &= w(t,1) = 0, \label{eq:target_bc}
\end{align}
where $K, c > 0$ are arbitrarily prescribed decay rates.

To analyze the exponential stability of \eqref{eq:target_ode}--\eqref{eq:target_bc}, we define the state space $\mathcal{H} \coloneqq \mathbb{C}\times L^{2}(0,1)$, endowed with the inner-product norm
\begin{align}
    \|(\vartheta,w)\|^{2}_{\mathcal{H}} = |\vartheta|^{2} + \int_{0}^{1} |w(x)|^{2}\,dx. \label{eq:norm_schrodinger}  
\end{align}

The dynamics of \eqref{eq:target_ode}--\eqref{eq:target_bc} can then be written in the abstract form
\begin{align}
    \frac{d Y_{w}(t)}{dt} &= \mathcal{A}_{w} Y_{w}(t), 
    \qquad 
    Y_{w}(0) = Y_{w0},
    \label{eq:abstract_evolution}
\end{align}
where $Y_{w}(t) = (\vartheta(t), w(t,\cdot)) \in \mathcal{H}$ and the operator 
$\mathcal{A}_{w}\!: D(\mathcal{A}_{w}) \subset \mathcal{H} \to \mathcal{H}$ 
is defined by
\begin{align}
    \mathcal{A}_{w}(\vartheta,w) = \big(-K\vartheta + w(0),\, -i w'' - c w\big),
    \label{eq:operator_schrodinger}
\end{align}
with domain
\begin{align}
    D(\mathcal{A}_{w}) =& 
    \Big\{(\vartheta,w)\in \mathbb{C}\times H^{2}(0,1)\ \big|\nonumber \\
    &\qquad \qquad \ w(1)=0,\ w'(0)=0 \Big\}.
    \label{eq:domain_operator_schrodinger}
\end{align}

This formulation allows the coupled ODE--PDE system to be viewed as an abstract evolution equation in a Hilbert space, which provides a convenient framework for the subsequent spectral and semigroup analysis. 

\medskip

\begin{lemma}
Let $\mathcal{A}_{w}$ be defined by \eqref{eq:operator_schrodinger}--\eqref{eq:domain_operator_schrodinger}. Then, the inverse $\mathcal{A}_{w}^{-1}$ exists and is compact on $\mathcal{H}$. Consequently, the spectrum of $\mathcal{A}_{w}$ consists only of isolated eigenvalues with finite algebraic multiplicity, explicitly given by 
\begin{align*} 
\lambda_{0} &= -K, & \lambda_{m} &= -c + i\left(m+\tfrac{1}{2}\right)^{2}\pi^{2}, \quad m\in\mathbb{N},
\end{align*}
    \noindent with eigenfunctions
\begin{align*}
    Z_0&=(\vartheta_0,0),\\
    Z_m(x)&=\big(\vartheta_m,\,w_m(x)\big),\quad m\in\mathbb{C},
\end{align*}
\noindent where $\vartheta_0,\,\vartheta_m\neq 0$ and $w_m(x)=\cos\!\big((m+\tfrac12)\pi x\big)$.
\end{lemma}

\medskip

\begin{pf}
\noindent\textbf{Existence and compactness of $\mathcal{A}_{w}^{-1}$.}  
Let $(\vartheta_{1},w_{1})\in\mathcal{H}$. The equation $\mathcal{A}_{w}(\vartheta,w)=(\vartheta_{1},w_{1})$ leads to the system
\begin{align}
    -K\vartheta+w(0) &= \vartheta_{1}, \label{eq:operator_schrodinger1}\\
    -iw''(x)-cw(x) &= w_{1}(x), \label{eq:operator_schrodinger2}\\
    w(1)=w'(0)&=0. \label{eq:operator_schrodinger3}
\end{align}

Equation \eqref{eq:operator_schrodinger2} is a second-order linear ODE with boundary conditions \eqref{eq:operator_schrodinger3}. By the variation of parameters formula, its unique solution is
\begin{align}
w(x) = \varphi(x) - \frac{\varphi(1)}{\cosh(\sqrt{c i})}\cosh(\sqrt{c i}\,x),
    \label{eq:inverse_operator_schrodinger1}
\end{align}
\noindent where 
\begin{align*}
    \varphi(x) = \sqrt{c i}\!\int_{0}^{x}\sinh\!\big(\sqrt{c i}(x-\tau)\big) w_{1}(\tau)\, d\tau.
\end{align*}

Once $w$ is obtained, $\vartheta$ follows from \eqref{eq:operator_schrodinger1}:
\begin{align}
    \vartheta = K^{-1}\big(w(0)-\vartheta_{1}\big). 
    \label{eq:inverse_operator_schrodinger2}
\end{align}

Hence, for each $(\vartheta_{1},w_{1})\in\mathcal{H}$ there exists a unique $(\vartheta,w)\in D(\mathcal{A}_{w})$.  
The mapping $\mathcal{A}_{w}^{-1}\!:\mathcal{H}\to D(\mathcal{A}_{w})$ is bounded, and by the Sobolev embedding theorem, the inclusion $D(\mathcal{A}_{w})\hookrightarrow\mathcal{H}$ is compact.  
Therefore, $\mathcal{A}_{w}^{-1}$ is compact on $\mathcal{H}$, and $\sigma(\mathcal{A}_{w})$ consists of isolated eigenvalues with finite algebraic multiplicity. 

\medskip

\noindent\textbf{Eigenvalues and eigenfunctions.}  
Let $\mathcal{A}_{w}Y_{w}=\lambda Y_{w}$ with $Y_{w}=(\vartheta,w)$. Then
\begin{align}
    -K\vartheta+w(0) &= \lambda \vartheta, \label{eq:eigen1}\\
    -iw''-cw &= \lambda w, \label{eq:eigen2}\\
    w(1)&=0,\quad w'(0)=0. \label{eq:eigen3}
\end{align}

Problem \eqref{eq:eigen2}-\eqref{eq:eigen3} has the simple spectrum
\begin{align*}
    \lambda_m=-c+i\mu_m,\qquad \mu_m=\Big(m+\tfrac12\Big)^2\pi^2,\quad m\in\mathbb{N},
\end{align*}
\noindent with eigenfunctions
\begin{align*}
   w_m(x)=\cos\!\Big(\big(m+\tfrac12\big)\pi x\Big).
\end{align*}

Substituting $\lambda_m$ and $w_m(0)=1$ into \eqref{eq:eigen1}, we obtain the corresponding eigenvector $\vartheta_m = 1/(K+\lambda_m)$. So $\lambda_m$, $m\in\mathbb{N}$, is the eigenvalue of \eqref{eq:eigen1}-\eqref{eq:eigen3} withe the corresponding eigenfunction $(1/(K+\lambda_m),\, w_m (x))$.

Next, let $\lambda = \lambda_0$. Since $-K\vartheta_{0}=\lambda_0 \vartheta_0$, it follows that $-K\vartheta_{0}+w(0)=\lambda_0 \vartheta_0$ yields $w(0)\equiv 0$. Moreover,
\begin{align*}
    -iw''-cw &= \lambda_0 w,\\
    w(0)&=w(1)= w'(0)=0.
\end{align*}
\noindent only has the trivial solution. So $\lambda_0$ is the eigenvalue of \eqref{eq:eigen1}-\eqref{eq:eigen3} and has the eigenfunction $(\vartheta_0 , \, 0)$. 
\end{pf}

\medskip

\begin{lemma}\label{prop:gen-stab}
Let $\mathcal A_w$ be given by \eqref{eq:operator_schrodinger}--\eqref{eq:domain_operator_schrodinger}. Then, there is a sequence of eigenfunctions of $\mathcal A_w$ which forms a Riesz basis for $\mathcal{H}$. Moreover,  $\mathcal A_w$ generates an exponentially stable $C_0$-semigroup on $\mathcal H$, that is, there exist $M\ge1$ and $\alpha\ge c$ such that
\begin{align*}
    \|e^{\mathcal A_w t}\|_{\mathcal L(\mathcal H)} \le M e^{-\alpha t},\qquad t\ge0.
\end{align*}
\end{lemma}

\medskip

\begin{pf}
Note that $\vartheta_0$ is an orthogonal basis for $\mathbb{C}$, and $\{w_{m}(x)=\cos((m+\tfrac{1}{2})\pi x),\ m\in\mathbb{N}\}$ is an orthonormal basis of $L^{2}(0,1)$. Then, $$\{w_0,\, w_m(x),\, m\in \mathbb{N}\},$$
where $F_0 = (\vartheta_0,\, 0)$ and $F_m(x) = (0,\, w_m (x))$, form a orthogonal basis in $\mathcal{H}$. Now, note that
\begin{align*}
    &|Z_0-F_0|^2 + \sum_{m=0}^{\infty}\|Z_{m}(x)-F_{m}\|^2 = \sum_{m=0}^{\infty}|(\lambda_{m}+K)|^{-2}\\
    & ~\quad\qquad\qquad\qquad\qquad\qquad\leq \frac{1}{\pi^4}\sum_{m=0}^{\infty}\frac{1}{(m+1/2)^4} = \frac{1}{6}.
\end{align*}
Therefore, by Bari's theorem, $\{Z_{0},\,Z_{m}(x)|\, m\in \mathbb{N}\}$ forms a Riesz basis for $\mathcal{H}$. Hence, there exists a bounded invertible $T:\ell^2\to\mathcal H$ such that $\mathcal A_w=T\Lambda T^{-1}$ with $$\Lambda=\mathrm{diag}(-K,\,-c+i\mu_0,\,-c+i\mu_1,\dots).$$ 
Since $\Lambda$ generates the diagonal semigroup $e^{\Lambda t}$ on $\ell^2$, we have
\begin{align*}
e^{\mathcal A_w t}=T\,e^{\Lambda t}\,T^{-1}\quad\Rightarrow\quad
\|e^{\mathcal A_w t}\|\le \|T\|\,\|T^{-1}\|\,\|e^{\Lambda t}\|.
\end{align*}

Because $\Re(-K)=-K$ and $\Re(-c+i\mu_m)=-c$ for all $m$, it follows that
\begin{align*}
    \|e^{\Lambda t}\|=\sup\{e^{-Kt},e^{-ct}\}=e^{-\min\{K,c\}\,t}.
\end{align*}

Therefore,
\begin{align*}
    \|e^{\mathcal A_w t}\|\le M_1\,e^{-\alpha t},
\end{align*}
\noindent where $M_1=\|T\|\,\|T^{-1}\|$ and $\alpha=\min\{K,c\}$.
\end{pf}

\medskip

\subsubsection{Backstepping transformation}

As highlighted in \cite{Ren2013}, applying a single-step backstepping transformation is challenging due to the complexity of the associated kernels. To overcome this difficulty, a two-step design approach is adopted.

\textbf{First backstepping transformation:} We consider the ODE-PDE system \eqref{eq:ODE1errordynamics}--\eqref{eq:BCUerrordynamics} and apply the backstepping transformation
\begin{align}
    v(t,x) = \beta(t,x) - \int_{0}^{x} q(x,y) \beta(t,y) \, dy - \gamma(x) \vartheta(t), \label{eq:backsteppingtransformation1}
\end{align}
\noindent which maps the original system \eqref{eq:ODE1errordynamics}--\eqref{eq:BCUerrordynamics} into the following intermediate target system:
\begin{align}
    \dot{\vartheta}(t) &= -K \vartheta(t) + v(t,0), \label{eq:intermediateODE} \\
    v_{t}(t,x) +i v_{xx}(t,x)&= 0, \label{eq:intermediatePDE} \\
    v_{x}(t,0) &= 0, \quad v(t,1) = W(t), \label{eq:intermediateBC1}
\end{align}
\noindent where $W \in \mathbb{C}$.

The kernels $q$ and $\gamma$ satisfy
\begin{align}
    \gamma''(x) &= 0, \label{eq:odekernel_intermediate} \\
    \gamma'(0) &= 0, \\ 
    \gamma(0) &= -K, \label{eq:odekernelbc2_intermediate} \\
    q_{yy}(x,y) - q_{xx}(x,y) &= 0, \label{eq:pdekernel_intermediate} \\
    q_{y}(x,0) &= -i \gamma(x), \quad q(x,x) = 0. \label{eq:pdekernelbc2_intermediate}
\end{align}

It is straightforward to verify that the solution of \eqref{eq:odekernel_intermediate}--\eqref{eq:odekernelbc2_intermediate} is
\begin{align}
    \gamma(x) = -K, \quad \forall x \in [0,1].
\end{align}

Moreover, the solution of \eqref{eq:pdekernel_intermediate}--\eqref{eq:pdekernelbc2_intermediate} is
\begin{align}
    q(x,y) = i \int_{0}^{x-y} \gamma(\sigma) \, d\sigma = -i K (x-y). \label{eq:solutionkernelpde_intermediate}
\end{align}

Finally, using \eqref{eq:backsteppingtransformation1} and \eqref{eq:intermediateBC1}, the resulting control law is
\begin{align}
    U(t) = W(t) - i K \int_{0}^{1} (1-y) \beta(t,y) \, dy - K \vartheta(t). \label{eq:controller_schrodinger}
\end{align}

\textbf{Second backstepping transformation:} Now, consider the following backstepping transformation:
\begin{align}
    w(t,x) = v(t,x) - \int_{0}^{x}\kappa(x,y)v(t,y)dy, \label{eq:backsteppingtransformation2}
\end{align}
\noindent to transform \eqref{eq:intermediateODE}-\eqref{eq:intermediateBC1} into \eqref{eq:target_ode}-\eqref{eq:target_bc}.

Differentiating \eqref{eq:backsteppingtransformation2} once with respect to time and twice with respect to space, substituting \eqref{eq:intermediateODE}-\eqref{eq:intermediateBC1} into it, and plugging the expressions into \eqref{eq:target_ode}-\eqref{eq:target_bc}, we obtain that \eqref{eq:intermediateODE}-\eqref{eq:intermediateBC1} is mapped into \eqref{eq:target_ode}-\eqref{eq:target_bc} if, and only if, the kernel $\kappa$ satisfies the following PDE:
\begin{align}
    &\kappa_{xx}(x,y)-\kappa_{yy}(x,y) = ic\kappa (x,y),\label{eq:kernelpdefinal}\\
    &\kappa_{y}(x,0) = 0, \quad
    \kappa(x,x) = -i\frac{c}{2}x.\label{eq:kernelbc2final}
\end{align}

The solution to the PDE \eqref{eq:kernelpdefinal}-\eqref{eq:kernelbc2final} is given in page 66 of \cite{Krsti2008BoundaryCO}, by
\begin{align}
    \kappa(x,y) = \kappa_{r}(x,y) + i\kappa_{i}(x,y),
\end{align}
\noindent where 
\begin{align*}
    \kappa_{r}(x,y) &= x\sqrt{\frac{c}{2(x^2-y^2)}}\left[-\ber_{1}\biggl(\sqrt{c(x^2-y^2)}\biggl) \right.\\
    & - \left. \bei_{1}\biggl(\sqrt{c(x^2-y^2)}\biggl)\right],\\
    \kappa_{i}(x,y) &= x\sqrt{\frac{c}{2(x^2-y^2)}}\left[\ber_{1}\biggl(\sqrt{c(x^2-y^2)}\biggl) \right.\\
    & - \left. \bei_{1}\biggl(\sqrt{c(x^2-y^2)}\biggl) \right],
\end{align*}
and $\ber_{1}$ and $\bei_{1}$ are the Kelvin functions.

Then, evaluating \eqref{eq:backsteppingtransformation2} at $x=1$ and using \eqref{eq:target_bc}, we obtain 
\begin{align}
    W(t) = \int_{0}^{1}\kappa(1,y)v(t,y)dy. \label{eq:intermediatecontroller}
\end{align}

Plugging \eqref{eq:intermediatecontroller} into \eqref{eq:controller_schrodinger}, with the help of \eqref{eq:backsteppingtransformation1}, we obtain the backstepping feedback control
\begin{align}
    U(t) &= \int_{0}^{1}\mathcal{S}_{1} (y)v(t,y)dy +K\mathcal{S}_{2}\vartheta(t), \label{eq:control_law_mean}
\end{align}
\noindent where
\begin{align*}
    \mathcal{S}_{1} (y) &= \left[\kappa_{r}+\int_{y}^{1}\kappa_{i}(1,\xi)(\xi-y)d\xi\right] \\
    &+i\biggl[\kappa_{i}(1,y)-K(1-y)
    -\int_{y}^{1}\kappa_{r}(1,\xi)(\xi-y)d\xi\biggr],\\
    \mathcal{S}_{2} &= \int_{0}^{1}\kappa(1,y)dy-1.
\end{align*}

\subsubsection{Invertibility of the Transformations}
As demonstrated in \cite{Ren2013}, the transformations \eqref{eq:backsteppingtransformation1} and \eqref{eq:backsteppingtransformation2} are invertible. Specifically, by postulating the inverse transformation of \eqref{eq:backsteppingtransformation1} as
\begin{align*}
    \beta(t,x) = v(t,x) - \int_{0}^{x} \xi(x,y) v(t,y) \, dy - \rho(x) \vartheta(t),
\end{align*}
\noindent and similarly, for \eqref{eq:backsteppingtransformation2},
\begin{align*}
    v(t,x) = w(t,x) - \int_{0}^{x} \eta(x,y) w(t,y) \, dy - \chi(x) \vartheta(t),
\end{align*}
one can determine the kernels $\xi$, $\eta$, $\rho$, and $\chi$ using the same reasoning as in the direct transformation. As a result, the closed-loop system and the target system exhibit the same stability properties.

\subsubsection{Implementable extremum seeking control law}

Note that the control law \eqref{eq:control_law_mean} cannot be directly applied in practice because the measurement of $\vartheta$ is not available. To circumvent this problem, we use the averaged version of the gradient and Hessian estimate (see \cite{book2022} for more details)
\begin{align}
    G_{av}(t)=H\vartheta_{av}(t), \qquad \hat{H}_{av}=H.\label{eq:gradient_estimateav}
\end{align}

In addition, recalling that $\psi$ denotes the isomorphism mapping a complex number to its real-imaginary representation (see Section \ref{section:control_problem}), we have that  $$z_{av}(t) = \frac{1}{\Pi}\int_{0}^{\Pi}\Gamma(s) M(s)y(s)ds = \Gamma_{av}(t)H\psi(\vartheta_{av}(t)),$$ where $\Pi = 2\pi \, \mbox{LCM}\{1/\omega_{i}\}$, for $i\in \{1,2\}$, and $\Gamma_{av}$ and $\vartheta_{av}$ denote the average version of $\Gamma$ and $\vartheta$, respectively. 

Since $\tilde{\Gamma}_{av}(t)=\Gamma_{av}(t)-H^{-1}$, it follows that
\begin{align}
     z_{av}(t) = \psi(\vartheta_{av}(t)) +\tilde{\Gamma}_{av}(t)H\psi(\vartheta_{av}(t)).
\end{align}

The linearization of the above expression at $H^{-1}$ results in $z_{av}(t)=\psi(\vartheta_{av}(t))$. Averaging \eqref{eq:control_law_mean} and using the above results, we obtain
\begin{align}
    U_{av}(t) &= \int_{0}^{1}\mathcal{S}_{1\,}(y)v_{av}(t,y)dy+K\, \vartheta_{av}(t). \label{eq:control_law_av}
\end{align}

For the stability analysis in which the averaging theory, for infinite-dimensional systems, is used, we employ a low-pass filter for the above feedback controller and then derive an infinite-dimensional averaging-based feedback given by 
\begin{align}
    U(t) = \mathcal{T}\left\{\int_{0}^{1}\mathcal{S}_{1}(y)v(t,y)dy+K\, \psi^{-1}(z(t)) \right\},\label{eq:filtered_controller}
\end{align}
\noindent where $\psi^{-1}$ is the inverse of $\psi$ (see Section \ref{section:control_problem}), and $\mathcal{T}$ is the low pass operator defined by
\begin{align*}
    \mathcal{T}\{g (t)\} = \mathcal{L}^{-1}\left\{ \frac{c_f}{s+c_f}\right\}*g (t),
\end{align*}
\noindent being $c_f>0$ the cutoff frequency, $\mathcal{L}^{-1}$ is the inverse Laplace transform, and $*$ is the convolution operation.

\section{Stability and convergence analysis}\label{section:stability}

The proof of stability is now presented. 

\medskip

\begin{theorem}
    For sufficiently large $c_f>0$, there exists $\bar{\omega}(c_f)>0$ such that $\forall \omega >\bar{\omega}$, the error dynamics has a unique exponentially stable periodic solution in $t$ of period $\Pi=2\pi\, \mbox{LCM}\{1/\omega_{i}\}$, for $i \in \{1,\,2\}$, denoted by $\tilde{\Gamma}^\Pi$, $\vartheta^\Pi$ and $\beta^\Pi(x,t)$, satisfying $\Psi(t) \leq \mathcal{O}(1/\omega)$, where
\begin{multline}
    \Psi(t) = \left(|\tilde{\Gamma}^\Pi(t)|^{2} + |\vartheta^\Pi(t)|^{2} + \int_{0}^{1}|\beta^\Pi(t,x)|^{2}dx\right.\\
    \left.+\int_{0}^{1}|\beta^\Pi_{x}(t,x)|^{2}dx\right)^{1/2}.\label{eq:psi-norm}
\end{multline}

Furthermore,
\begin{align}
    \lim\sup_{t\rightarrow\infty}|\theta(t)-\Theta^{*}| &= \mathcal{O}\left(|a|\mathrm{e}^{\sqrt{\omega}}+ 1/\omega\right),\label{eq:convergence_theta}\\
    \lim\sup_{t\rightarrow\infty}|\Theta(t)-\Theta^{*}| &= \mathcal{O}(|a| +1/\omega),\label{eq:convergence_Theta}\\
    \lim\sup_{t\rightarrow\infty}|y(t)-y^{*}| &= \mathcal{O}(|a|^{2}+1/\omega^{2}),\label{eq:convergence_y}
\end{align}
\noindent where $|a|=\sqrt{a_{1}^{2}+a_{2}^{2}}$.
\end{theorem}
\begin{pf} 
\textbf{Average Model of the Closed-loop System.} Recalling that $\psi^{-1}(z_{av}(t))=\vartheta_{av}(t)$, we obtain, for $\omega$ sufficiently large, the averaged version of the closed-loop system~ \eqref{eq:ODE1errordynamics}–\eqref{eq:BCUerrordynamics}, with control law \eqref{eq:filtered_controller}:
\begin{align}
    &\dot{\vartheta}_{av}(t) = \beta_{av}(t,0), \label{eq:ODE1errordynamics_av}\\
    &(\beta_{av})_{t}(t,x) + i(\beta_{av})_{xx}(t,x) = 0, \label{eq:PDEerrordynamics_av}\\
    &(\beta_{av})_{x}(t,0) = 0, \label{eq:BC0errordynamics_av}\\
    &\dot{\beta}_{av}(t,1) = - c_{f}\beta_{av}(t,1) +
    c_f\int_0^1\mathcal{S}_1(y)v_{av}(t,y)dy \nonumber\\
 &\qquad\qquad\!\! +c_f\mathcal{S}_2 \vartheta_{av}(t),\label{eq:BCUerrordynamics_av}\\
 &v_{av}(t,x) = \beta_{av}(t,x) + i\int_{0}^{x}K(x-y)\beta_{av}(t,y)dy\nonumber\\
 &\qquad\qquad\!\! +K\vartheta_{av}(t).\label{eq:control_law_av}
\end{align}

\textbf{Backstepping transformation and target system.} Using backstepping transformations analogous to \eqref{eq:backsteppingtransformation1}–\eqref{eq:backsteppingtransformation2}, the averaged dynamics \eqref{eq:ODE1errordynamics_av}–\eqref{eq:control_law_av} can be mapped accordingly:
\begin{align}
    &\dot{\vartheta}_{av}(t) = -K\vartheta_{av}(t) + w_{av}(t,0), \label{eq:target_ode_av}\\
    &(w_{av})_{t}(t,x) + i (w_{av})_{xx}(t,x)= - c w_{av}(t,x), \label{eq:target_pde_av}\\
    &(w_{av})_{x}(t,0) = 0, \label{eq:target_bc_av}\\
    &w_{av}(t,1) = - c_{f}w_{av}(t,1) +
    c_f\int_0^1\mathcal{S}_1(y)v_{av}(t,y)dy \nonumber\\
 &\qquad\qquad\!\! +c_f\mathcal{S}_2\vartheta_{av}(t),\label{eq:BCUerrordynamics_av}\\
 &v_{av}(t,x) = w_{av}(t,x) + i\int_{0}^{x}K(x-y)w_{av}(t,y)dy\nonumber\\
 &\qquad\qquad\!\! +K\vartheta_{av}(t).
\end{align}
For clarity, we denote the averaged transformed states by $w_{av}$ and $v_{av}$, distinguishing them from the previously defined variables $w$ and $v$.

\textbf{Lyapunov stability of the target system.} Consider the Lyapunov functional
\begin{align}
    V(t) &= |\vartheta_{av}(t)|^{2} + \frac{p_{1}}{2}\!\int_{0}^{1}\!|w_{av}(t,x)|^{2}dx \nonumber\\
    & + \frac{p_{2}}{2}\!\int_{0}^{1}\!|(w_{av})_{x}(t,x)|^{2}dx,
    \label{eq:lyap_schrodinger}
\end{align}
where $p_{1},p_{2}>0$ are constants to be determined.

Differentiating \eqref{eq:lyap_schrodinger} along the trajectories of 
\eqref{eq:target_ode}–\eqref{eq:target_bc} yields, after assuming $c_f\to\infty$ in \eqref{eq:filtered_controller} for simplicity,
\begin{align}
    \dot{V}(t) &= -2K|\vartheta_{av}(t)|^{2}
    + 2\Re\{w_{av}(t,0)\overline{\vartheta_{av}(t)}\}
    \nonumber\\
    &- c\!\int_{0}^{1}\!(p_{1}|w_{av}(t,x)|^{2}+p_{2}|(w_{av})_{x}(t,x)|^{2})dx \nonumber\\
    & + \int_{0}^{1}\!\!p_{1}\Im\{(w_{av})_{xx}(t,x)\overline{w_{av}(t,x)}\}dx\nonumber\\
    & +\int_{0}^{1}\!\! p_{2}\Im\{(w_{av})_{xxx}(t,x)\overline{(w_{av})_{x}(t,x)}\}dx.
    \label{eq:diff_lyap_schrodinger}
\end{align}

By integrating by parts and using the boundary conditions 
$(w_{av})_{x}(t,0)=w_{av}(t,1)=0$, one verifies that the last two integrals vanish. Hence
\begin{align*}
    \dot{V}(t)
    &= -2K|\vartheta_{av}(t)|^{2}
    + 2\Re\{w_{av}(t,0)\overline{\vartheta_{av}(t)}\}\\
    &- c\!\int_{0}^{1}\!(p_{1}|w_{av}(t,x)|^{2}+p_{2}|(w_{av})_{x}(t,x)|^{2})dx. 
\end{align*}

Then, applying Young's and Agmon's inequalities, it follows that there exists $\gamma > 0$ such that
\begin{align*}
    \dot{V}(t) \leq& \left(\frac{1}{\gamma}-K\right)|\vartheta_{av} (t)|^{2}+\left(\gamma-p_{2}\,c\right)\int_{0}^{1}|w_{av}(t,x)|^{2}dx\\
    & + \left(\gamma-p_{3}\,c\right)\int_{0}^{1}|(w_{av})_{x}(t,x)|^{2}dx.
\end{align*}
    Then, choosing $p_{1},p_{2}>\frac{1+c^2}{c^{2}}$ and $\alpha = \min \{K, c\}$, we obtain $\dot{V}\leq -\alpha V(t),$ or
    \begin{align}
        V(t)\leq V(0)\mathrm{e}^{-\alpha t}.\label{eq:lyap_inequality}
    \end{align}



Moreover, the linearized average model of the Hessian inverse estimation error around $\tilde{\Gamma}_{av} = 0$ satisfies $\dot{\tilde{\Gamma}}_{av}(t) = -\omega_r \tilde{\Gamma}_{av}(t)$, which is exponentially stable for any $\omega_r > 0$. Consequently, using these results and the invertibility of the backstepping transformations~\eqref{eq:backsteppingtransformation1}--\eqref{eq:backsteppingtransformation2}, the averaged closed-loop system~\eqref{eq:lift_cl_pde}--\eqref{eq:vcompact} is exponentially stable.

\textbf{Invoking averaging theorem.} Let $v(t,x):=\beta(t,x)-U(t)$. The coupled ODE–PDE dynamics \eqref{eq:ODE1errordynamics}–\eqref{eq:BCUerrordynamics} can be written as
\begin{align}
\dot{\vartheta}&=v(t,0)+U(t),\label{eq:lift_cl_pde}\\
v_t(t,x)+i\,v_{xx}(t,x)&=-\phi(U,v),\\
v_x(t,0)=v(t,1)&=0,\\ 
\dot{U}(t)&=\phi(U,v),\\
\dot{\tilde{\Gamma}}(t)&=\omega_r\tilde{\Gamma}(t)-\omega_r\Gamma(t)\hat{H}(t)\tilde{\Gamma}(t), \label{eq:vcompact}
\end{align}
where
\begin{align*}
\phi(U,v)&=-c_fU+c_f\!\int_0^1\!\mathcal{S}_1(y)\!\left(v(t,y)+U(t)\right)\!dy\\
&+c_f\mathcal{S}_2 \psi^{-1}(z(t)).
\end{align*}

Let $\alpha(v) = \int_0^1S_1(y)v(y)dy$, $s_1 = \int_0^1S_1(y)dy$, and define the linear operator
\begin{align*}
    \mathcal{A}_{v}(\vartheta,U,\Gamma,v) = (v(0)+U, \, - c_f U,\, \omega_r\Gamma,\, iv''+c_fU),
\end{align*}
\noindent with domain
\begin{align*}
    D(\mathcal{A}_{v}) &= 
    \Big\{(\vartheta,U,\Gamma,v)\in \mathbb{C}^{2}\times\mathbb{R}^{2}\times H^{2}(0,1)\ \big|\nonumber \\
    &\qquad \qquad \qquad \qquad \ v(1)=0,\ v'(0)=0 \Big\}.
\end{align*}

Then, introducing 
$Y_v(t)=(\vartheta(t),\;U(t),\;\Gamma(t), \;v(t,\cdot))$, the closed-loop system \eqref{eq:lift_cl_pde}-\eqref{eq:vcompact} takes the abstract form
\begin{equation}
\frac{dY_v}{dt}(t)=\mathcal{A}_v Y_v(t)+J(\omega t,Y_v(t)), \label{eq:abstractcompact}
\end{equation}
with
\begin{align*}
J(\omega t,Y_v)=\left(\begin{array}{c}
 0,\\
c_f\left(\alpha(v) + s_1 U(t) + \mathcal{S}_2G(t)\right)\\
-\omega_r\Gamma(t)\hat{H}(t)\Gamma(t)\\
-c_f\left(\alpha(v) + s_1 U(t) + \mathcal{S}_2G(t)\right)
\end{array}\right).
\end{align*}

The operator $\mathcal{A}_v$ can be written as a matrix operator composed of the infinite-dimensional block $\mathcal{A}_w$ and the finite-dimensional dynamics of $U$ and $\Gamma$. Since $\mathcal{A}_w$ and $\mathrm{diag}(-c_f,\omega_r)$ each generate a $C_0$-semigroup on their respective spaces, the composite operator $\mathcal{A}_v$ generates a $C_0$-semigroup on $\mathbb{C}^2\times\mathbb{R}^2\times L^2(0,1)$. As the PDE block has compact resolvent, $\mathcal{A}_v$ also has compact resolvent, so its spectrum is discrete with finite algebraic multiplicity. Moreover, the nonlinear term $J(\omega t,Y_v)$ is Fréchet differentiable in $Y_v$, strongly continuous, and almost periodic in $t$, fulfilling the assumptions of the infinite-dimensional averaging theorem~\cite{book2022}. Hence, the closed-loop system~\eqref{eq:lift_cl_pde}--\eqref{eq:vcompact} admits a periodic solution. Then, using \eqref{eq:lyap_inequality} and the fact that $\tilde{\Gamma}_{av}$ is exponentially stable, it follows that, for $\omega$ sufficiently large, the closed-loop system ~\eqref{eq:lift_cl_pde}--\eqref{eq:vcompact}, with $U$ in \eqref{eq:filtered_controller}, has a unique exponentially stable periodic solution around its equilibrium
satisfying $\Psi(t)\leq \mathcal{O}(1/\omega)$, with $\Psi$ defined in \eqref{eq:psi-norm}.

\textbf{Asymptotic convergence to a neighborhood of the
extremum.} The asymptotic convergence to a neighborhood of the extremum point follows from 
$|\theta(t) - \Theta^{*}| = |\tilde{\theta}(t) + S(t)|$, with $S$ in \eqref{eq:perturbation} being of order $\mathcal{O}(|a|\mathrm{e}^{\sqrt{\omega}})$, leading directly to 
\eqref{eq:convergence_theta}. Following the procedure in Theorem~10.1 (step~6) of~\cite{book2022}, one can analogously derive \eqref{eq:convergence_Theta}. Finally, from \eqref{extrachapter.eq:static_map}, we have $|y(t)-y^{*}|\leq \lambda_{\max}(H)|\Theta(t)-\Theta^{*}|^{2}$, where $\lambda_{\max}(\cdot)$ denotes the largest eigenvalue. Taking the limit superior, substituting \eqref{eq:convergence_Theta}, and applying Young’s inequality yields \eqref{eq:convergence_y}.
\end{pf}

\section{Numerical Results}\label{section:results}

\begin{figure}[t!]
    \centering
    \includegraphics[width=\linewidth]{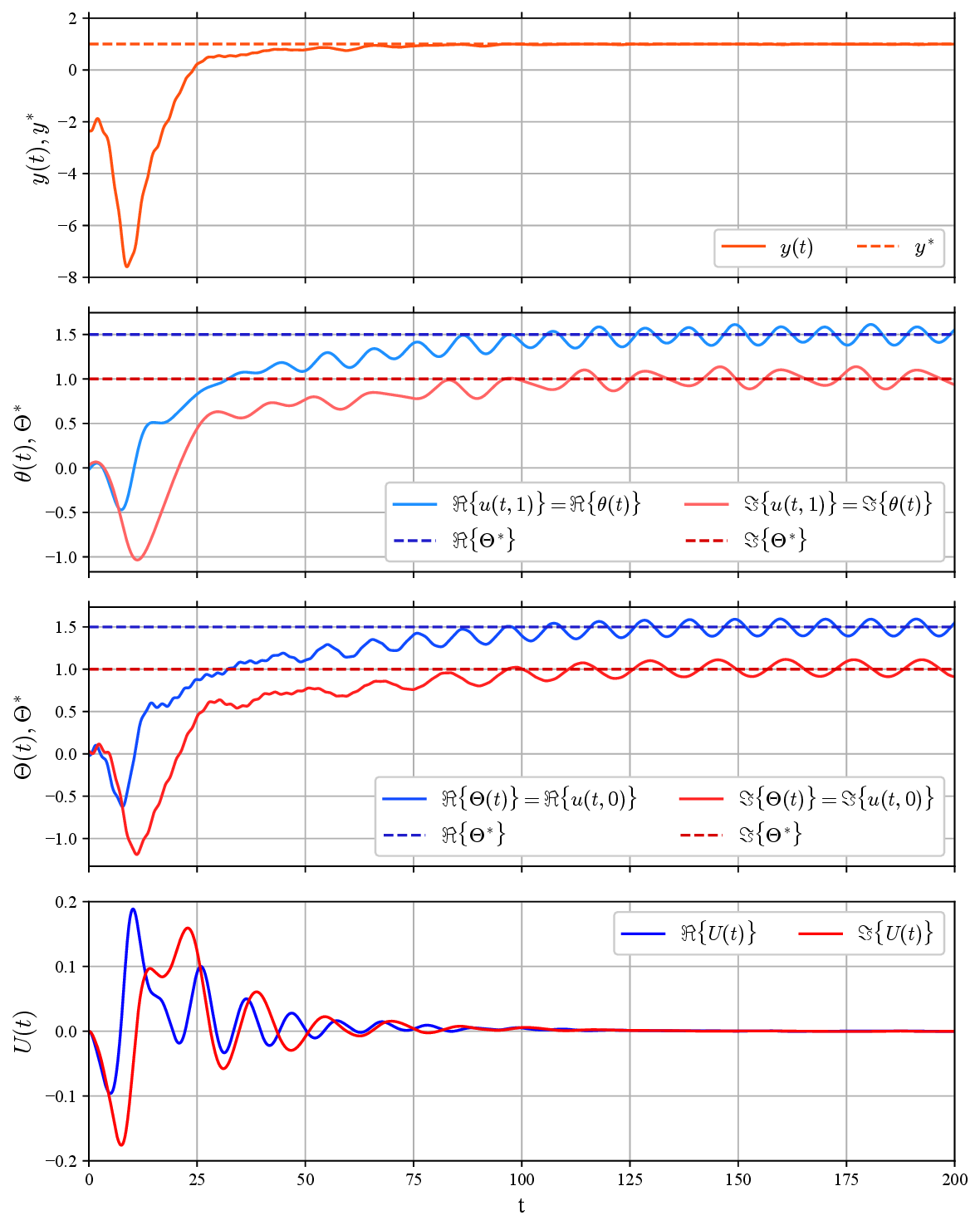}
    \caption{The closed-loop response of the Schrödinger PDE with ES compensating controller.}
    \label{fig:results}
\end{figure}

\begin{figure}[t]
    \centering
    \includegraphics[width=0.95\linewidth]{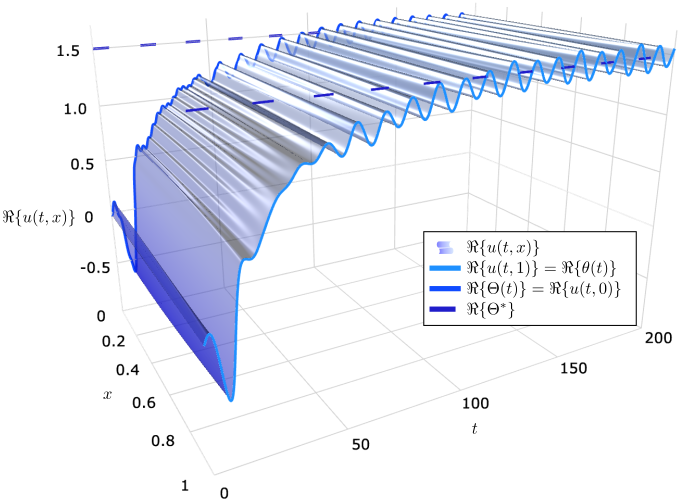}
    \caption{The response of the real part of the Schrödinger equation evolving in time and space.}
    \label{fig:3Dr}
\end{figure}

\begin{figure}[t]
    \centering
    \includegraphics[width=0.95\linewidth]{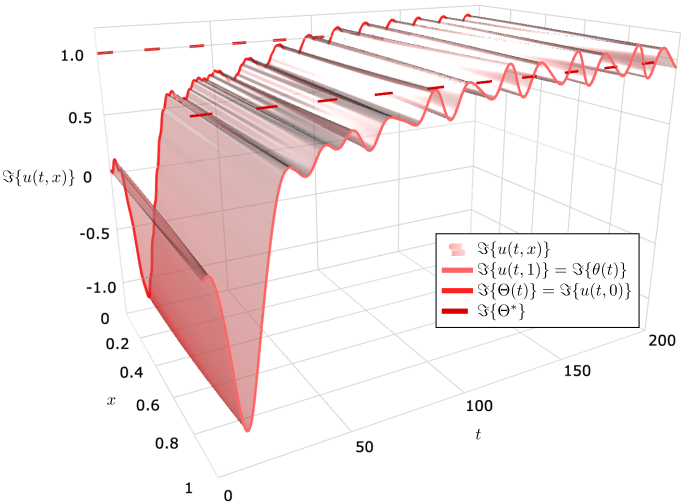}
    \caption{The response of the imaginary part of the Schrödinger equation evolving in time and space.}
    \label{fig:3Di}
\end{figure}

The numerical implementation of the linearized Schrödinger equation was carried out using the finite element method with Hermitian functions. Considering a quadratic static map as in \eqref{extrachapter.eq:static_map}, the system is subjected to the control law \eqref{eq:control_law_av}. The Hessian is given by $H=-\begin{bmatrix} 1.0&1.0\\1.0&1.5 \end{bmatrix}$, with an optimizer $\Theta^{*} = 1.5+i$ and an optimal unknown output value $y^{*}=1.0$. The controller parameters are chosen as $\omega_1 = 0.6$, $\omega_2 = 0.4$, $\omega_r = 0.7$, $a_1=a_2=0.1$, $c=0.1$, $c_f = 50$, and $K = 0.13$. The closed-loop simulation results, presented in Figure \ref{fig:results}, demonstrate the effectiveness of the proposed control approach. The control actions depicted in Figure \ref{fig:results} ensure that the variables $(y,\theta,\Theta)$ converge toward the neighborhood of their optimal values $(y^{*},\Theta^{*},\Theta^{*})$. These findings validate the performance of the ES-based control strategy in driving the system towards optimal operation as shown in Figures \ref{fig:3Dr} and \ref{fig:3Di}.

\section{Conclusions}\label{section:conclusion}

The proposed extremum seeking framework optimizes a quadratic performance map linked to the boundary value of a linearized Schrödinger equation. This identifies the optimal steady state without prior knowledge of map parameters, using only real-time measurements. A boundary feedback law based on backstepping compensates for the infinite-dimensional actuation dynamics, enabling average-based estimation of the map’s gradient, Hessian and its inverse. The closed-loop system converges to a neighborhood of the optimum while maintaining exponential stability of the Schrödinger dynamics, thus bridging data-driven adaptive optimization and boundary control of dispersive PDEs using extremum seeking.

\bibliographystyle{plain}
\bibliography{bibli}

\end{document}